\documentclass[12pt]{article}
\usepackage{amsmath}
\usepackage{amssymb}

\usepackage[cp1251]{inputenc}
\usepackage[russian]{babel}
\newcommand{\il}[2]{\int\limits_{#1}^{#2}}

\newcommand{\ph}{\phantom{a}}
\newcommand{\phh}{\phantom{aaa}}

\newcommand{\sist}[2]{\left\{
\begin{array}{l}
{#1}\\
\ph\\
{#2}
\end{array}
\right.}

\begin{document}
\vskip 12pt\noindent
MSC 34K99, 34K43
\vskip 12pt
\centerline {\bf G. A. Grigorian}
\vskip 12pt
{\bf\centerline {Cauchy problem for quasilinear systems}
\centerline {of functional differential equations}}

\vskip 10 pt

\centerline{0019 Armenia c. Yerevan, str. M. Bagramian 24/5}
\centerline{Institute of Mathematics of NAS of Armenia}
\centerline{E - mail: mathphys2@instmath.sci.am, \ph phone: 098 62 03 05, \ph 010 35 48 61}

\vskip 20 pt

\noindent
Abstract. We use the contracting mapping principle for proving that under some mild restrictions the Cauchy problem for quasilinear systens of functional differential equations with retarded arguments has the unique solution. As a consequence from this result we obtain that the Cauchy problem for linear systems of functional differential equations with locally integrable coefficients and with locally measurable retarded arguments has the unique solution. We show that similar results can be obtained for the Cauchy co problem of quasilinear systems of functional differential equations with advanced arguments.

\vskip 12pt
\noindent
Key words: quasilinear systems of functional differential equations, the Cauchy problem, the Cauchy co problem, retarded arguments, advanced arguments, the contracting mapping principle, linear systems of functional differential equations.

\vskip 12pt

{\bf 1. Introduction}. Let  $F_k(t,u_{11},\dots,u_{1N},\dots,u_{nN}), \ph k=\overline{1,n}$ be real-valued locally integrable by $t$ and continuous by $u_{11},\dots,u_{nN}$ functions on $[t_0,+\infty)\times \mathbb{R}^{nN}$, and let $\alpha_{kj}(t), \ph k=\overline{1,n}, \ph j=\overline{1,N}$ be real-valued locally measurable functions on $[t_0,+\infty)$.
 Consider the  system of functional differential equations
$$
\phi'(t) = F_k(t,\phi_1(\alpha_{11}(t)),\dots,\phi_1(\alpha_{1N}(t)),\dots,\phi_n(\alpha_{n1}(t)),\dots,\phi_n(\alpha_{nN}(t))), \ph t \ge t_0,  \eqno (1.1)
$$
$ k = \overline{1,n}$. Let $r_k(t), \ph k = \overline{1,n}$ be real-valued continuous functions on $(-\infty,t_0]$. By a Cauchy problem for the system (1.1) we mean to find a real-valued continuous vector function $(\phi_1(t),\dots,\phi_n(t))$ on $\mathbb{R}$, which is absolutely continuous on $[t_0,+\infty)$, and which satisfies (1.1) almost everywhere on $[t_0,+\infty)$ and the initial conditions
$$
\phi_k(t) = r_k(t), \ph t \le t_0, \ph k=\overline{1,n}. \eqno (1.2)
$$
Throughout of this paper we will assume that the following conditions are satisfied:

\noindent
(L) \ph (the Lipshits's condition)
$$
|F_k(t,u_{11},\dots,u_{nN}) - F_k(t,v_{11},\dots,v_{nN})| \le f_k(t)\sum\limits_{m=1}^n\sum\limits_{j=1}^N |u_{mj} - v_{mj}|, \ph t \ge t_0, \ph u_{mj},  v_{mj} \in \mathbb{R},
$$
where $f_k(t), \ph k=\overline{1,n}$ are locally integrable functions on $[t_0,+\infty)$;

\noindent
(LI) (local integrability)

for any locally measurable functions $\psi_{11}(t),\dots,\psi_{1N}(t),\dots,\psi_{n1}(t),\dots,\psi_{nN}(t)$  on $[t_0, + \infty)$ the  superpositions
$$
F_k(t,\psi_{11}(t),\dots,\psi_{nN}(t)),\phh k=\overline{1,n}
$$
are locally integrable on $[t_0,+\infty)$;

\noindent
(Ret) (the retorsion conditions)
$$
\alpha_{kj}(t) \le t, \phh t \ge t_0, \phh k = \overline{1,n}, \phh j = \overline{1,N}.
$$

{\bf Remark 1.1.} {\it The condition (LI) is satisfied if in particular
$$
F_k(t,u_{11},\dots,u_{nN}) = \sum\limits_{m=1}^Mf_{km}(t)g_{km}(u_{11},\dots,u_{nN}), \ph t \ge t_0, \ph u_{11},\dots,u_{nN} \in \mathbb{R}
$$
where $f_{km}(t), \ph k=\overline{1,n}, \ph m=\overline{1,M}$ are locally integrable functions on $[t_0,+\infty)$, \linebreak $g_{km}(u_{11},\dots,u_{nN}), \ph \ph k=\overline{1,n}, \ph m=\overline{1,M}$ are continuous functions on $\mathbb{R}^{nN}$.}

Let  $G_k(t,u_{11},\dots,u_{1N},\dots,u_{nN}), \ph k=\overline{1,n}$ be real-valued locally integrable by $t$ and continuous by $u_{11},\dots,u_{nN}$ functions on $(-\infty, \tau_0]\times \mathbb{R}^{nN}$, and let $\beta_{kj}(t), \ph k=\overline{1,n}, \ph j=\overline{1,N}$ be real-valued locally measurable functions on $(-\infty, \tau_0]$.
 Consider the system of functional differential equations
$$
\phi'(t) = G_k(t,\phi_1(\beta_{11}(t)),\dots,\phi_1(\beta_{1N}(t)),\dots,\phi_n(\beta_{n1}(t)),\dots,\phi_n(\beta_{nN}(t))), \ph t \le \tau_0,  \eqno (1.3)
$$
$ k = \overline{1,n}$. Let $s_k(t), \ph k = \overline{1,n}$ be real-valued continuous functions on $(-\infty,t_0]$. By a Cauchy co problem for the system (1.3) we mean to find a real-valued continuous vector function $(\phi_1(t),\dots,\phi_n(t))$ on $\mathbb{R}$, which is absolutely continuous on $(-\infty, \tau_0]$, and which satisfies (1.3) almost everywhere on $(-\infty, \tau_0]$  and the terminal conditions
$$
\phi_k(t) = s_k(t), \ph t \ge \tau_0, \ph k=\overline{1,n}. \eqno (1.4)
$$
Throughout of this paper we will assume that the following conditions are satisfied:

\noindent
$(L^*)$ \ph (the Lipshits's condition)
$$
|G_k(t,u_{11},\dots,u_{nN}) - G_k(t,v_{11},\dots,v_{nN})| \le h_k(t)\sum\limits_{m=1}^n\sum\limits_{j=1}^N |u_{mj} - v_{mj}|, \ph t \ge t_0, \ph u_{mj},  v_{mj} \in \mathbb{R},
$$
where $h_k(t), \ph k=\overline{1,n}$ are locally integrable functions on $(-\infty, \tau_0]$;

\noindent
($LI^*$) (local integrability)

for any locally measurable functions $\psi_{11}(t),\dots,\psi_{1N}(t),\dots,\psi_{n1}(t),\dots,\psi_{nN}(t)$  on $(-\infty, \tau_0]$  the \linebreak superpositions
$$
G_k(t,\psi_{11}(t),\dots,\psi_{nN}(t)),\phh k=\overline{1,n}
$$
are locally integrable on $(-\infty, \tau_0]$;

\noindent
(Adv) (the advance   conditions)
$$
\beta_{kj}(t) \ge t, \phh t \le \tau_0, \phh k = \overline{1,n}, \phh j = \overline{1,N}.
$$

Functional differential equations and systems of such equations appear in various areas of natural science, such as in Economics (see, e.g. [1-6]), in the probability theory (see, e.g., [7-9]), in Biology (see, e.g., the predator-prey model of Volterra [11, p. 3], the model of  circumnutation of plants [11, p.3],  the model of dynamics of individuals, infected by gonorrhea  [11, p. 4]), in the electrodynamics [11, p. 7], so on. Therefore the qualitative study of functional differential equations  is very actual.  For qualitative study the solutions of these equations and systems of such equations the main interest represents the case when the studying solution is continuable on whole semi axis.
Therefore, the study the Cauchy problem (1.1), (1.2) (the Cauchy co problem (1.3), (1.4)) is  very actual as well. Probably the Cauchy problem is studied systematically and solved for linear systems of functional differential equations with coefficients from wide classes of functions but with specified deviations of the argument of the form $t + \xi_k, \ph  \xi_k = const, \ph k=1,2,\dots$ (see, e.g., [11, 12]) and for some particular classes of equations with deviations of the argument from wide classes of functions (see, e.g., [12-14]). However, probably, in most of cases the global existence of solutions of studying equations and (or) systems of equations is only assumed, making this way their results  conditional (see, e.g., [15-25]).

In this paper we show that the Cauchy problem (1.1), (1.2) (as well as the cauchy co problem (1.3), (1.4)) under some mild restrictions has the unique solution, which weakens (to some extent) the conditionality of some results, in which the global existence of solutions of studying equations and systems of equations is only assumed.

\vskip 12pt

{\bf 2. Auxiliary propositions}. Let $t_0 \le t_1 < t_2 < +\infty$ and let $\psi_1(t),\dots,\psi_n(t)$ be real-valued continuous functions on $(-\infty,t_2]$. Denote by $AC_{\psi_1,\dots,\psi_n}^n[t_1,t_2]$ the set of all real-valued continuous vector functions $x(t)\equiv (x_1(t),\dots,x_n(t))$ on $(-\infty,t_2]$ such that $x(t)$ is absolutely continuous on $[t_1,t_2]$ and $x_k(t) = \psi_k(t), \ph t \le t_1, \ph k=\overline{1,n}.$ Obviously, $AC_{\psi_1,\dots,\psi_n}^n[t_1,t_2]$ is a closed (full) metric space with the distance $\rho(x,t) \equiv \sum\limits_{k=1}^n\max\limits_{t\in[t_1,t_2]}|x_k(t) - y_k(t)|$ between its elements $x(t) \equiv (x_1(t),\dots,x_n(t)), \ph y(t)\equiv (y_1(t),\dots,y_n(t)).$ For any $x(t)\equiv (x_1(t),\dots,x_n(t)) \in AC_{\psi_1,\dots,\psi_n}^n[t_1,t_2]$ set
$$
(I_{\psi_1,\dots,\psi_n,t_1,t_2} x)(t) \equiv \sist{x(t), \phh t \le t_1,}{((I_{1,\psi_1,\dots,\psi_n,t_1,t_2}x)(t),\dots,(I_{n,\psi_1,\dots,\psi_n,t_1,t_2}x)(t)), \ph t\in [t_1,t_2],}
$$
where $(I_{k,\psi_1,\dots,\psi_n,t_1,t_2}x)(t)\equiv \psi_k(t_1) +$

\phantom{a} $+ \il{t_1}{t}F_k(\tau,x_1(\alpha_{11}(\tau)),\dots,x_1(\alpha_{1N}(\tau)),
\dots,x_n(\alpha_{n1}(\tau)),\dots,x_n(\alpha_{nN}(\tau)))d\tau,\ph t\in [t_1,t_2],$ $k=\overline{1,n}.$
Obviously under the restrictions (LI) and (Ret) the operator $I_{\psi_1,\dots,\psi_n,t_1,t_2}$ is a mapping from $AC_{\psi_1,\dots,\psi_n}^n[t_1,t_2]$ into itself.

{\bf Theorem 2.1.} {\it Let the conditions (L), (LI) and (Ret) be satisfied. If \linebreak $N\sum\limits_{k=1}^{n}\il{t_1}{t_2}f_k(\tau)d\tau <~1$, then $I_{\psi_1,\dots,\psi_n,t_1,t_2}$ is a contracting mapping in $AC_{\psi_1,\dots,\psi_n}^n[t_1,t_2]$.}

Proof. It follows from (L) and (LI) that for every $x(t)\equiv (x_1(t),\dots,x_n(t)), \ph y(t)\equiv (y_1(t),\dots,y_n(t)) \in AC_{\psi_1,\dots,\psi_n}^n[t_1,t_2]$ the following chin of relations is valid.
$$
\rho(I_{\psi_1,\dots,\psi_n,t_1,t_2} x, I_{\psi_1,\dots,\psi_n,t_1,t_2} y) = \sum\limits_{k=1}^{n}\max\limits_{t\in[t_1,t_2]}\biggl|\il{t_1}{t}F_k(\tau,x_1(\alpha_{11}(\tau)),\dots,x_n(\alpha_{nN}(\tau))) d\tau - \phantom{aaaaaa}
$$
$$
\phantom{aaaaaaaaaaaaaaaaaaaaaaaaaaaaaaaaaa}- \il{t_1}{t}F_k(\tau,y_1(\alpha_{11}(\tau)),\dots,y_n(\alpha_{nN}(\tau))) d\tau\biggr | \le
$$
$$
\sum\limits_{k=1}^{n}\max\limits_{t\in[t_1,t_2]}\il{t_1}{t}f_k(\tau)\biggr\{\sum\limits_{m=1}^n
\sum\limits_{j=1}^N|x_m(\alpha_{mj}(\tau)) - y_m(\alpha_{mj}(\tau))|\biggr\} d \tau \le  \phantom{aaaaaaaaaaaaaaaaaaaaaaaaa}
$$
$$
\phantom{aaaaaaaaaaaaaaaaaaaaaaaaaaaaaaaaaaaaaa}\le \sum\limits_{k=1}^{n}\il{t_1}{t_2}f_k(\tau) d\tau D(x,y), \eqno (2.1)
$$
where $D(x,y) \equiv \max\limits_{\tau\in[t_1,t_2]}\sum\limits_{m=1}^n\sum\limits_{j=1}^N|x_m(\alpha_{mj}(\tau)) - y_m(\alpha_{mj}(\tau))|.$ From (Ret) it follows that $D(x,y) \le N \rho(x,y).$ This together with (2.1) implies that $\rho(I_{\psi_1,\dots,\psi_n,t_1,t_2} x, I_{\psi_1,\dots,\psi_n,t_1,t_2} y) \le \linebreak \le N\sum\limits_{k=1}^{n}\il{t_1}{t_2}f_k(\tau)d\tau \rho(x,y).$ Therefore, if $N\sum\limits_{k=1}^{n}\il{t_1}{t_2}f_k(\tau)d\tau <~1$, then $I_{\psi_1,\dots,\psi_n,t_1,t_2}$ is a contracting mapping in $AC_{\psi_1,\dots,\psi_n}^n[t_1,t_2]$. The theorem is proved.

 Let $-\infty <  t_1 < t_2 \le \tau_0$ and let $\chi_1(t),\dots,\chi_n(t)$ be real-valued continuous functions on $[t_2,+\infty)$. Denote by $AC_{\chi_1,\dots,\chi_n}^{n*}[t_1,t_2]$ the set of all real-valued continuous vector functions $x(t)\equiv (x_1(t),\dots,x_n(t))$ on $[t_1,+\infty)$ such that $x(t)$ is absolutely continuous on $[t_1,t_2]$ and $x_k(t) = \chi_k(t), \ph t \ge t_2, \ph k=\overline{1,n}.$ Obviously, $AC_{\chi_1,\dots,\chi_n}^{n*}[t_1,t_2]$ is a closed (full) metric space with the distance $\rho(x,t) \equiv \sum\limits_{k=1}^n\max\limits_{t\in[t_1,t_2]}|x_k(t) - y_k(t)|$ between its elements $x(t) \equiv (x_1(t),\dots,x_n(t)), \ph y(t)\equiv (y_1(t),\dots,y_n(t)).$ For any $x(t)\equiv (x_1(t),\dots,x_n(t)) \in \linebreak AC_{\psi_1,\dots,\psi_n}^{n*}[t_1,t_2]$ set
$$
(J_{\chi_1,\dots,\chi_n,t_1,t_2} x)(t) \equiv \sist{x(t), \phh t \ge t_2,}{((J_{1,\chi_1,\dots,\chi_n,t_1,t_2}x)(t),\dots,(J_{n,\chi_1,\dots,\chi_n,t_1,t_2}x)(t)), \ph t\in [t_1,t_2],}
$$
where $(J_{k,\chi_1,\dots,\chi_n,t_1,t_2}x)(t)\equiv \chi_k(t_2) +$

\phantom{a} $+ \il{t}{t_2}G_k(\tau,x_1(\beta_{11}(\tau)),\dots,x_1(\beta_{1N}(\tau)),
\dots,x_n(\beta_{n1}(\tau)),\dots,x_n(\beta_{nN}(\tau)))d\tau,\ph t\in [t_1,t_2],$ $k=\overline{1,n}.$
Obviously under the restrictions ($LI^*$) and (Adv) the operator $J_{\chi_1,\dots,\chi_n,t_1,t_2}$ is a mapping from $AC_{\psi_1,\dots,\psi_n}^{n*}[t_1,t_2]$ into itself. By analogy of the proof of Theorem 2.1 it can be proved

{\bf Theorem 2.2.} {\it Let the conditions ($L^*$), ($LI^*$) and (Adv) be satisfied. If \linebreak $N\sum\limits_{k=1}^{n}\il{t_1}{t_2}h_k(\tau)d\tau <~1$, then $J_{\chi_1,\dots,\chi_n,t_1,t_2}$ is a contracting mapping in $AC_{\psi_1,\dots,\psi_n}^{n*}[t_1,t_2]$.}

\vskip 12pt

{\bf 3. Main results.} By a solution of the system (1.1) on $[t_1,t_2)\ph (t_0\le t_1< t_2 \le \infty)$ we mean a real-valued continuous vector function $(\phi_1(t),\dots,\phi_n(t))$ on $(-\infty,t_2)$, which is absolutely continuous on $[t_1,t_2)$ and satisfies (1.1) almost everywhere on $[t_1,t_2)$.

{\bf Definition 3,1,} {\it An interval $[t_0,T)\ph (t_0 < T \le +\infty)$ is called the maximum existence interval for a solution $\Phi(t)$ of the system (1.1), if $\Phi(t)$ exists on $[t_0,T)$ and cannot be continued to the right from $T$ as a solution of (1.1).}

{\bf Theorem 3.1.} {\it Let the conditions (L), (LI) and (Ret) be satisfied. Then the Cauchy problem (1.1), (1.2) has the unique solution.}

Proof. Since  $f_k(t), \ph k=\overline{1,n}$ are locally integrable chose $t_1 > t_0$ so close to $t_0$ that
$$
N\sum\limits_{k=1}^n\il{t_0}{t_1}f_k(\tau) d\tau <1.
$$
Then by Theorem 2.1 it follows from (L), (LI) and (Ret) that the operator $I_{r_1,\dots,r_n,t_0,t_1}$ is a contracting mapping in $AC_{r_1,\dots,r_n}^n[t_0,t_1]$. Therefore, according to the contracting mapping principle, the operator $I_{r_1,\dots,r_n,t_0,t_1}$ has the unique fixed point $\Phi(t)\equiv (\phi_1(t),\dots,\phi_n(t))$ in $AC_{r_1,\dots,r_n}^n[t_0,t_1]$, which is a solution of the system (1.1) on $[t_0,t_1)$, satisfying the initial conditions $\phi_k(t) = r_k(t), \ph t\le t_0, \ph k= \overline{1,n}$. It follows from here that $\Phi(t)$ has a maximum existence interval. Let $[t_0,T)$ be that interval. Show that
$$
T = +\infty. \eqno (3.1)
$$
Suppose
$$
T < +\infty. \eqno (3.2)
$$
Chose $\varepsilon > 0$ so small that
$$
N\sum\limits_{k=1}^n\il{T-\varepsilon}{T + \varepsilon}f_k(\tau) d\tau <1.
$$
Then by Theorem 2.1  it follows from (L), (LI) and (Ret) that  $I_{\phi_1,\dots,\phi_n,T-\varepsilon,T+\varepsilon}$ is a contracting mapping in $AC_{r_1,\dots,r_n}^n[T-\varepsilon,T+\varepsilon]$, and, according to the contracting mapping principle, has the unique fixed point $\widetilde{\Phi}(t)\equiv (\widetilde{\phi_1}(t),\dots,\widetilde{\phi_n}(t))$ in  $AC_{r_1,\dots,r_n}^n[T-\varepsilon,T+\varepsilon]$, which is a solution of the system (1.1) on $[T-\varepsilon,T+\varepsilon)$, satisfying the initial conditions
$$
\widetilde{\phi_k}(t) = \phi_k(t), \ph t \le T - \varepsilon, \ph k= \overline{1,n}.
$$
This means (since $\widetilde{\Phi}(t)$ is the unique) that $\widetilde{\Phi}(t)$ coincides with $\Phi(t)$ on $(-\infty,T]$ and is a continuation of $\Phi(t)$ on $(-\infty,T+\varepsilon]$ as a solution of the system (1.1) on $[t_0,T+\varepsilon)$, which contradicts (3.2). The obtained contradiction proves (3.1). So, to complete the proof of the theorem it remains to show the uniqueness of $\Phi(t)$. Suppose there exists another, different from $\Phi(t)$ solution $\Psi(t)\equiv (\psi_1(t),\dots,\psi_n(t))$ to the problem (1.1), (1.2). Then there exist $t_0\le T_1 < T_2 < +\infty$ such that
$$
\Phi(t) = \Psi(t), \ph t \le T_1, \ph \Phi(t) \ne \Psi(t), \ph t \in (T_1,t_2). \eqno (3.3)
$$
Without loss of generality we may take that
$$
N\sum\limits_{k=1}^n\il{T_1}{T_2}f_k(\tau) d\tau <1.
$$
Then by virtue of Theorem 2.1 and the contracting mapping principle the operator $I_{\phi_1,\dots,\phi_n,T_1,T_2}$ has the unique fixed point $(v_1(t),\dots,v_n(t))$  in   $AC_{r_1,\dots,r_n}^n[T_1,T_2]$, which is the unique solution of the system (1.1) on $[T_1,T_2)$, satisfying the initial conditions
$$
v_k(t) = \phi_k(t), \ph t \le T_1, \ph k= \overline{1,n}.
$$
We obtain a contradiction with (3.3). The obtained contradiction completes the proof of the theorem.

Let $a_{kjm}(t)(t), \ph b_k(t), \ph k,j = \overline{1,n}, \ph m = \overline{1,N}$ be real-valued locally integrable functions on $[t_0,+\infty)$. Consider the linear system of functional differential equations
$$
\phi_k'(t) = \sum\limits_{j=1}^n\sum\limits_{m=1}^N a_{kjm}(t)\phi_j(\alpha_{jm}(t)) + b_k(t), \phh t \ge t_0, \ph k=\overline{1,n}. \eqno (3.4)
$$
This system is a particular case of (1.1), for which the conditions (L) and (LI) are satisfied. Then from Theorem 3.1 we immediately obtain

{\bf Corollary 3.1.} {\it Let (Ret) be satisfied. Then the Cauchy problem (3.4), (1.2) has the unique solution.}

Using Theorem 2.2 instead of Theorem 2.1 by analogy of the proof of Theorem 3.1 it can be proved

{\bf Theorem 3.2.} {\it Let the conditions ($L^*$), ($LI^*$) and (Adv) be satisfied. Then the Cauchy co problem (1.3), (1.4) has the unique solution.}

Let $c_{kjm}(t)(t), \ph d_k(t), \ph k,j = \overline{1,n}, \ph m = \overline{1,N}$ be real-valued locally integrable functions on $(-\infty,\tau_0]$. Consider the linear system of functional differential equations
$$
\phi_k'(t) = \sum\limits_{j=1}^n\sum\limits_{m=1}^N c_{kjm}(t)\phi_j(\beta_{jm}(t)) + d_k(t), \phh t \le \tau_0, \ph k=\overline{1,n}. \eqno (3.5)
$$
This system is a particular case of (1.3), for which the conditions   ($L^*$)and ($LI^*$) are satisfied. Then from Theorem 3.2 we immediately obtain

{\bf Corollary 3.1.} {\it Let (Adv) be satisfied. Then the Cauchy co problem (3.5), (1.4) has the unique solution.}

\vskip 30 pt

\centerline {\bf References}
\vskip 12pt

\noindent
1.R. Frisch and H. Holme, The Characteristic Solutions of a Mixed Difference and Differential \linebreak \phantom{a} Equations Occurring in Economic Dynamics. Economica, III (1935)  pp. 225 - 239.

\noindent
2. M. Kalecki, A Macrodinamic Theory of Business Cycles. Economica, III (1935)  \linebreak \phantom{a} pp.  327 - 344.

\noindent
3. A. Callender, D. R. Hertree, A. Porter, Time-Lag in a Control System, Philos. Trans.  \linebreak \phantom{a} Roy. Soc. London (A), 235 (766) (1936) pp. 415 - 444.

\noindent
4. D. R. Hartree, A. Porter, A. Callender, A. B. Stevenson, Time Lag in a Control System,  \linebreak \phantom{a} II. Proc. Roy. Soc. London (A) 161 (907) (1937) pp. 460 - 476.

\noindent
5. N. Minorsky, Control Problems, II, Jpurn. Frankl. Inst., 232: 6 (1941) pp. 519 - 551.

\noindent
6. H. Bateman, The Control of an Elastic Fluid. Bull. Amer. Math. Soc., 51 (1945)  \linebreak \phantom{a} pp. 601 - 641.

\noindent
7. L. Silberstein, On a Histero-Differential Equations, Arising in a Probability Problem.  \linebreak \phantom{a} Phil. Mag. (t), 29: 192 (1840) pp. 75 - 84.

\noindent
8. N. Hale, On the Statistical Treatment of Counting Experiments in Nuclear Physics.  \linebreak \phantom{a} Ark. Mat. Astr. Fys., 33A: 11 (1946) pp 1 - 11.

\noindent
9. N. Hole, On the Distribution on Counts in a Counting Apparatus. Ark. Mat. Astr. Fys., \linebreak \phantom{a}33:3 (1947) B:8, pp. 1 -8.

\noindent
10. N. Minorsky, Self - Exited Mechanical oscillations. Journ. Appl. Phys, 19 (1948) \linebreak \phantom{aa}  pp. 332 - 338.

\noindent
11. J. Hale, Theory of Functional Differential Equations. Applied Mathematical \linebreak \phantom{aa} Sciences | Vol. 3, Spriger-Verlag, New York, Heidelberg, Berlin, 1977, 366 pages.

\noindent
12. R. Bellman, K. Cooke, Differential-Difference Equations. Academic Press, New York,\linebreak \phantom{aa}  London, 1963, 465 pages.

\noindent
13. L. Berezansky and E. Braverman, Some Oscillation Problems for Second Order Linear \linebreak \phantom{aa} Delay Differential Equations. J. Math. Anal., Appl.., 220, 719-740 (1009)

\noindent
14. G. A. Grigorian,  Oscillation criteria for the second order linear functional-differential \linebreak \phantom{aa} equations with locally integrable coefficients. Saraevo Journal of Mathematics.  Vol.14 \linebreak \phantom{aa} (27), No.1, (2018), pp. 71–86.

\noindent
15. J. Dzurina, Oscillation theorems for second order advanced neutral differential equations. \linebreak \phantom{aa}  Tatra Mt. Math. Publ. 48 (2011) pp. 61 - 71.

\noindent
16. R. Guo, Q. Hung and Q. Liu, Some New Oscillation Criteria of Even-Order  \linebreak \phantom{aa} Quasi-linear Delay Differential Equations with Neutral Term, Mathematics, \linebreak \phantom{aa} 2021, 9, 2074, pp. 1-11.

\noindent
17. M. Pasic, Parametric Exited Oscillation of Second-Order Functional-Differential  \linebreak \phantom{aa} Equations and Application to Duffing Equations with Time Delay Feedback. Discrete  \linebreak \phantom{aa} Dynamics in Nature and Society, 2014, pp. 1-17.

\noindent
18. M. M. El-sheikh, R. Sallam and N. Mohamady. On the oscillation of third order neutral  \linebreak \phantom{aa} delay differential equations. Appl. Math. Inf. Sci. Lett., 1, No. 3, (2-13) pp. 77- 80.

\noindent
19. X. Lin, Oscillation of second-order neutral differential equations. J. Math. Anal.  Appl. \linebreak \phantom{aa} 309 (2005) pp. 442-452.

\noindent
20. H. - J. Li, Ch. - Ch. Yen, Oscillation of Nonlinear Functional Differential Equations  \linebreak \phantom{aa} of the second order. Appl. Math. Lett., vol. 11, No. 1, (1998) pp. 71-77.

\noindent
21. Sh. Tang, T. Li, E. Thamdapani, Oscillation of higher-order half-linear neutral  \linebreak \phantom{aa}  differential   equations. Demostrato Mathematica, Vol. XLVI, No. 1, 2013, pp. 101-109.

\noindent
22. A. A. Soliman, R. A. Sallam, A. M. Hassan, Oscillation Criteria of Seco0nd Order  \linebreak \phantom{aa} Nonlinear Neutral Differential Equations. Int. Journal of Applied Mathematical  \linebreak \phantom{aa} Research, 1 (3) (2012) pp. 314-322.

\noindent
23. J. R. Graef, S. Murugadas, E. Thandapani, Oscillation Criteria for Second Order  \linebreak \phantom{aa} Neutral Delay Differential Equations with Mixed Nonlinearities. International  \linebreak \phantom{aa} Electronic Journal of Pure and Applied Mathematics, Vol. 2. No. 1,  2010, pp. 85-99.

\noindent
24. P. Wang, Zh. Xu, On the oscillation of a two-dimensional delay differential systems.  \linebreak \phantom{aa} Int. Journ. of Qualitative Theory of Differential Equations and Applications, vol. 2., \linebreak \phantom{aa} No. 1 (2008) pp. 38-46.

\noindent
25. R. P.  Agraval, M. Bohner, T. Li, Ch. Zheng, Oscillation of second-order Emden-Fowler  \linebreak \phantom{aa} neutral delay differential equations. Annali di Matematica. 193 (2014) pp. 1861-1875.

\end{document}